\documentclass{article}
\usepackage[dvips]{epsfig}
\usepackage{color}

\usepackage[utf8]{inputenc}
\usepackage{amsmath, amssymb}
\usepackage[amsmath, thmmarks]{ntheorem}
\usepackage{graphicx}
\newtheorem{theorem}{Theorem}

\newtheorem{corollary}{Corollary}
\newtheorem{lemma}{Lemma}

\newenvironment{proof}{\par\noindent {\bf Proof}}{\hfill qed\par}

\newcommand{\bR}{\mathbb{R}}
\newcommand{\bN}{\mathbb{N}}
\def\dim {\ dim}
\def\e {\ e}
\def\pS{$\sf{S}$}
\def\SAT{${\sf 3-SATISFIABILITY}$}
\def\NAE3SAT{${\sf NOT~ALL~EQUAL~3-SATISFIABILITY}$}
\def\COL{${\sf 3-COLORABILITY}$}
\def\disjunion {\ \widetilde{\cup}\ }
\def\union {\ \cup\ }
\def\join {\ +\ }
\begin{document}

\title{On the computational complexity of degenerate unit-distance representations of graphs}

\author{
  Jan Kratochv\'il\thanks{\texttt{honza@kam.mff.cuni.cz}}\\
  Charles University, Prague, \\
  Czech Republic\\[12pt]
  Boris Horvat\thanks{\texttt{Boris.Horvat@fmf.uni-lj.si}}\\
  IMFM, University of Ljubljana, \\
  Slovenia\\[12pt]
  Toma\v{z} Pisanski\thanks{\texttt{Tomaz.Pisanski@fmf.uni-lj.si}} \\
  IMFM, University of Ljubljana, and University of Primorska, \\
  Slovenia
}
\date{\today}

\maketitle

%%%%%%%%%%%%%%%%%%%%%%%%%%%%%%
%
%  Abstract
%
%%%%%%%%%%%%%%%%%%%%%%%%%%%%%%

\begin{abstract}
Some graphs admit drawings in the Euclidean $k$-space in such a (natural) way, that edges are represented as line segments of unit length. Such drawings will be called $k$ dimensional unit distance representations. When two non-adjacent vertices are drawn in the same point, we say that the representation is degenerate. 
The dimension (the Euclidean dimension) of a graph is defined to be the minimum integer $k$ needed that a given graph has non-degenerate $k$ dimensional unit distance representation (with the property that non-adjacent vertices are mapped to points, that are not distance one appart). 

It is proved that deciding if an input graph is homomorphic to a graph with dimension $k \geq 2$ (with the Euclidean dimension $k \geq 2$) are NP-hard problems.\\ 

{\bf Mathematics Subject Classification}: 
05C62,  % Graph representations (geometric and intersection representations, etc.)
05C12,  % Distance in graphs
05C10,  % Topological graph theory, imbedding [See also 57M15, 57M25]
% http://www.ams.org/msc/index.html#search
\end{abstract}

%%%%%%%%%%%%%%%%%%%%%%%%%%%%%%
%
%  Introduction and background
%
%%%%%%%%%%%%%%%%%%%%%%%%%%%%%%

\section{Introduction and background}

% BH there are lot of different definitions that define the same thing in unit distance graph world. Different authors are using same terms for different things. Professor Pisanski and I decided to use terms such as unit distance representation, unit distance realization, unit distance graph and strict unit distance graph (or unit distance coordinatization). Those are the ``most popular'' or the most ``self-explanatory''.

% Representation, realization, degenerate
Define a {\em representation} $\rho := (\rho_{V},\rho_{E})$ of a graph $G$ in a set $M$ as a mapping $\rho_{V}: V(G) \ \rightarrow M$ and a mapping $\rho_{E}: E(G) \ \rightarrow 2^{M}$, such that if $v$ is an end-vertex of an edge $e=u \sim v$, then $\rho_{V}(v) \in \rho_{E}(e)$ (and $\rho_{V}(u) \in \rho_{E}(e)$) \cite{pearlsgt}. If the converse is also true, that is, that for $v \in V(G)$ and $e \in E(G)$, if $\rho_{V}(v) \in \rho_{E}(e)$ is true, then $v$ is an end-vertex of $e=u \sim v$, the representation is called a {\em realization}. If there is no danger of confusion we drop the subscripts and denote both mappings $\rho_{V}$ and $\rho_{E}$ by $\rho$. Observe, that the intersection $\rho(e) \cap \rho(\ell)$ of realizations of non-adjacent edges $e, \ell \in E(G)$ can be non-empty, but cannot contain a realization of any vertex of  $G$. Thus, representations of edges of $G$ may cross, which is different from the usual graph embedding case. If the mapping $\rho$ is not injective on $V(G)$ it is called {\em degenerate}. 

% Euclidean representation
In this paper we consider only representations in the Euclidean $k$-space, i.e. ({\em $k$ dimensional Euclidean}) representations that have $M = \bR^k$, and where an edge is always represented by the line segment between the representations of its end-vertices. Each Euclidean realization of a connected graph on at least three vertices is non-degenerate, see \cite{HorvatPisanski2}.
An Euclidean representation is called a {\em unit distance representation} if each edge is represented as a line segment of length one. A graph is called a $k$ dimensional {\em unit distance graph} if it has a $k$ dimensional unit distance realization. For example, the well known Petersen graph $G(5,2)$ is a unit distance graph with exactly 18 different degenerate unit distance representations, see \cite{HorvatPisanski}. 

% unit distance embedding, Dimension, Euclidean dimension
Erd\"os {\em et al.}, see \cite{ErdosHararyTutte}, defined the {\em dimension} of a graph $G$, denoted as $\dim(G)$, as the minimum integer $k$ needed that $G$ is a $k$ dimensional unit distance graph. Other authors (e.g. \cite{BuckleyHarary, UDGMAEHARA4, UDGMAEHARA2}) defined an {\em unit distance embedding} to be a unit distance realization that maps non-adjacent vertices to points with distances other than one. The smallest integer $k$ needed that a given graph $G$ has an unit distance embedding in $\bR^k$ is called the {\em Euclidean dimension} of $G$, see for example \cite{UDGMAEHARA3}, and is, as usual, denoted as $\e(G)$. A graph with $\e(G) = k$ is called a $k$ dimensional {\em strict unit distance graph} (or, following Boben {\em et al.} \cite{BobenPisanski}, $k$ dimensional {\em unit distance coordinatization}). As any $k$ dimensional strict unit distance graph is a $k$ dimensional unit distance graph, it follows that for a given graph $G$, $\dim(G) \leq \e(G)$. Wheh $\dim(G)=2$ (or $\e(G)=2$) we will be talking about {\em planar} unit distance graphs (or planar strict unit distance graphs). A planar graph with a planar unit distance realization in the plane is called a {\em matchstick} graph. 

% BH Now, when using terms dimension and Euclidean dimension, we are speaking with the terms that are commonly used. Euclidean embedding is the induced subgraph of the well known graph of the $k$ space. (Comes from the well known problem: how many colors are needed to color the plane?)

% Observing just simple graphs
If $G$ is composed of the connected components $G_{1}, G_{2}, \ldots, G_{s}$, then $G$ admits a unit distance representation in $\bR^{k}$ if and only if each component $G_{i}$, $1 \leq i \leq s$, admits a unit distance representation in $\bR^{k}$. By definition, loops are not valid unit distance Euclidean edge representations, and there exists only one line segment of length one between two points in the Euclidean space. Therefore, we may restrict our attention to connected simple graphs.

% homomorphism
A mapping $f: V(G) \rightarrow V(H)$ from a graph $G$ into a graph $H$ is called a {\em graph homorphism}, if $f$ maps vertices of $G$ into vertices of $H$, such that $u \sim v \in E(G)$ implies $f(u) \sim f(v) \in E(H)$ ($f$ preserves adjacencies). Each graph morphism $f: V(G) \rightarrow V(H)$ induces a unique mapping (denoted by the same letter) $f: E(G) \rightarrow E(H)$, such that for $e=u \sim v$, $f(e)=f(u)\sim f(v)$ (see e.g. \cite{GAHHELL}). Let $G$ and $H$ be graphs, let $f:G \rightarrow H$ be a graph homomorphism, and let $k \in \bN$. Any unit distance representation of $H$ in $\bR^{k}$ (if it exists) can be extended to a (possibly degenerate) unit distance representation of $G$ in $\bR^{k}$, see \cite{HorvatPisanski2}.

% graph join, union
The union $G_{1} \union G_{2}$ of graphs $G_{1}$ and $G_{2}$ is the graph with vertex set $V(G_{1}) \cup V(G_{2})$ and edge set  $E(G_{1}) \cup E(G_{2})$. The disjoint union of two graphs $G_{1}$ and $G_{2}$, denoted $G_{1} \disjunion G_{2}$, is the graph obtained by taking the union of $G_{1}$ and $G_{2}$ on disjoint vertex sets, $V(G_{1})$ and $V(G_{2})$. The join of two simple graphs $G_{1}$ and $G_{2}$, written $G_{1} \join G_{2}$  is the graph obtained by taking the disjoint union of $G_{1}$ and $G_{2}$ and adding all edges joining $V(G_{1})$ and $V(G_{2})$. 
%BH graph join naj bi bil zvezdica, ne plus?

\begin{figure}[ht]
\centerline{\epsfig{file=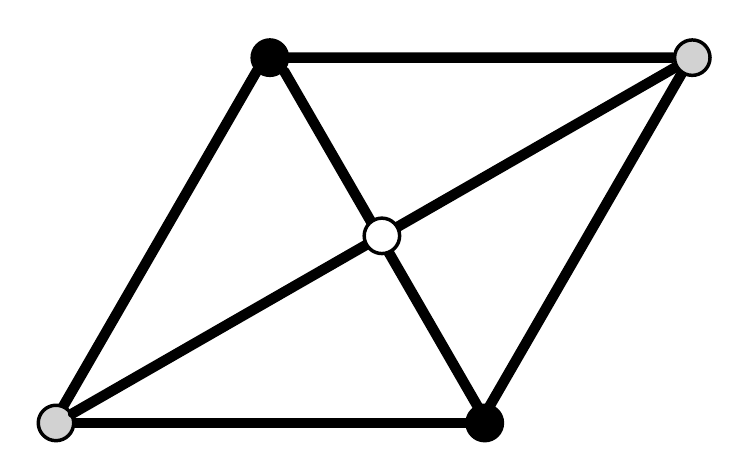,width=4cm} \epsfig{file=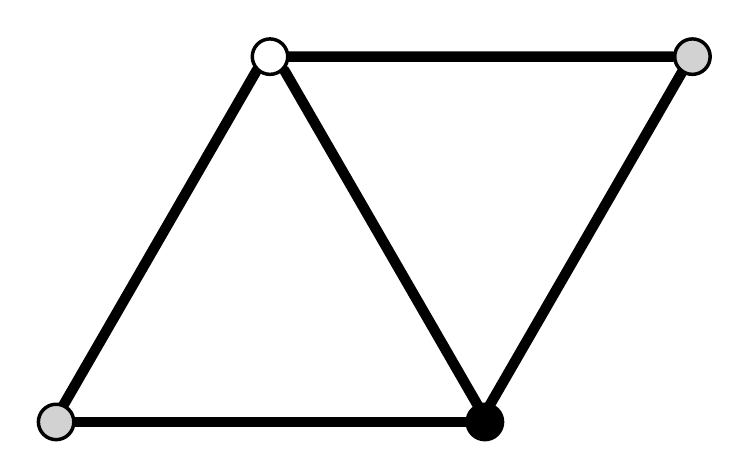,width=4cm}  \epsfig{file=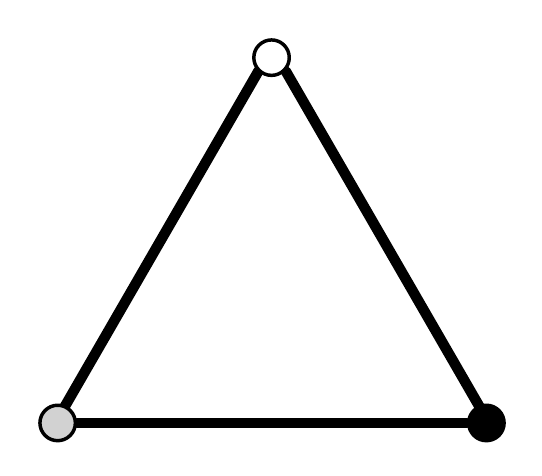,width=2.9cm} }
\caption{The wheel graph $W_{5}$ (on the left), which admits a proper 3 coloring, is not a planar unit distance graph. Degenerate planar unit distance representations of the wheel graph $W_{5}$ are obtained using black vertex identification (in the middle) and additional gray vertex identification (on the right).}
\label{fig:degudgrepr}
\end{figure}

A wheel graph $W_{n}$ on $n \geq 4$ vertices is defined as $W_{n} = W_{1, n-1} := K_{1} \join C_{n-1}$, where $K_{1}$ is the one vertex graf and $C_{n-1}$ the cycle graph on $n-1$ vertices. It is well known \cite{BuckleyHarary} that $\dim(W_{n}) = 3$ for $n \neq 7$ and $\dim (W_{7}) = 2$. 
Even though $\e(W_{5}) = \dim(W_{5}) = 3$, there exists degenerate planar unit distance representations of the graph $W_{5}$; see Figure~\ref{fig:degudgrepr}.

It is therefore interesting to answer a question ``Does, for a given graph $G$ and a natural number $k$, there exists a degenerate $k$ dimensional unit distance representation of $G$?''. One can instead, answer a question ``For a given natural number $k$ and a graph $G$, is $G$ homomorphic to a graph with dimension $k$?''. Of course, similar questions concerning the Euclidean dimension can be asked.

The only graph with zero dimension is $K_{1}$, and the same holds for the Euclidean dimension. It can be easily shown that the only graphs with dimension one are path graphs on at least two vertices. Again, the same is true for the Euclidean dimension. Thus, we can check in polynomial time whether a given graph is homomorphic to a linear strict unit distance graph. In Section 2 we prove that deciding if an input graph is homomorphic to a planar unit distance graph (planar strict unit distance graph) are NP-hard problems. In Section 3 we consider the same problems in higher dimensions.

%%%%%%%%%%%%%%%%%%%%%%%%%%%%%%
%
%  Planar unit distance representations
%
%%%%%%%%%%%%%%%%%%%%%%%%%%%%%%

\section{Planar unit distance representations}

% Theorem 1
\begin{theorem}
Deciding if an input graph $G$ is homomorphic to a graph $H$ with $\dim(H) = 2$ (with $\e(H) = 2$) are NP-hard problems.
\end{theorem}

\begin{figure}[ht]
\centerline{\epsfig{file=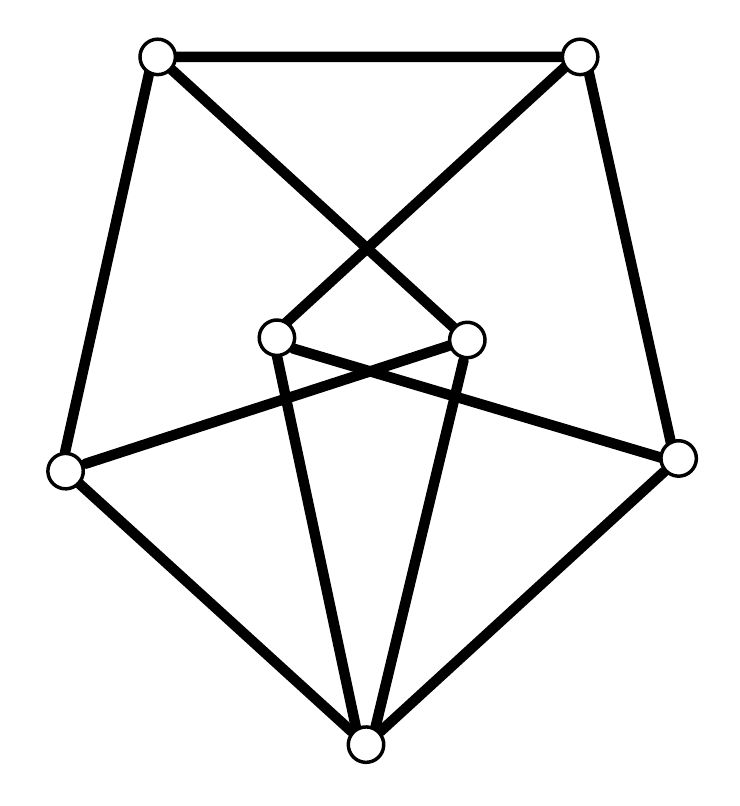,width=4cm}}
\caption{The well known Moser graph (Moser spindle) is a planar strict unit distance graph.}
\label{fig:Moser}
\end{figure}

\begin{proof}
The proof is based on the Moser graph (see Fig.~\ref{fig:Moser}) which is a well known planar strict unit distance graph. It has the property that every homomorphism to a unit distance graph $H$ is injective and the image $H$ is isomorphic to the Moser graph itself. In particular, adding any edge to the Moser graph results in a graph with dimension at least 3. It is not hard to see, that the dimension and the Euclidean dimension of the Moser graph are both equal to two. Using Laman's Theorem, see e.g. \cite{Hendrickson:1992p3225}, it is not hard to show that the unit distance coordinatization of Moser graph is rigid.

The NP-hardness reduction goes from \SAT. Given a formula $\Phi$ with a set of variables $X$ and a set of clauses $C$ (each clause containing exactly 3 literals), we construct a graph $G_{\Phi}$, such that $G_{\Phi}$ is homomomorphic to a strict unit distance graph (namely to the Moser graph) if $\Phi$ is satisfiable, while $G_{\Phi}$ is not homomorphic to any unit distance graph if $\Phi$ is unsatisfiable. This will prove both claims.

\begin{figure}[ht]
\centerline{\epsfig{file=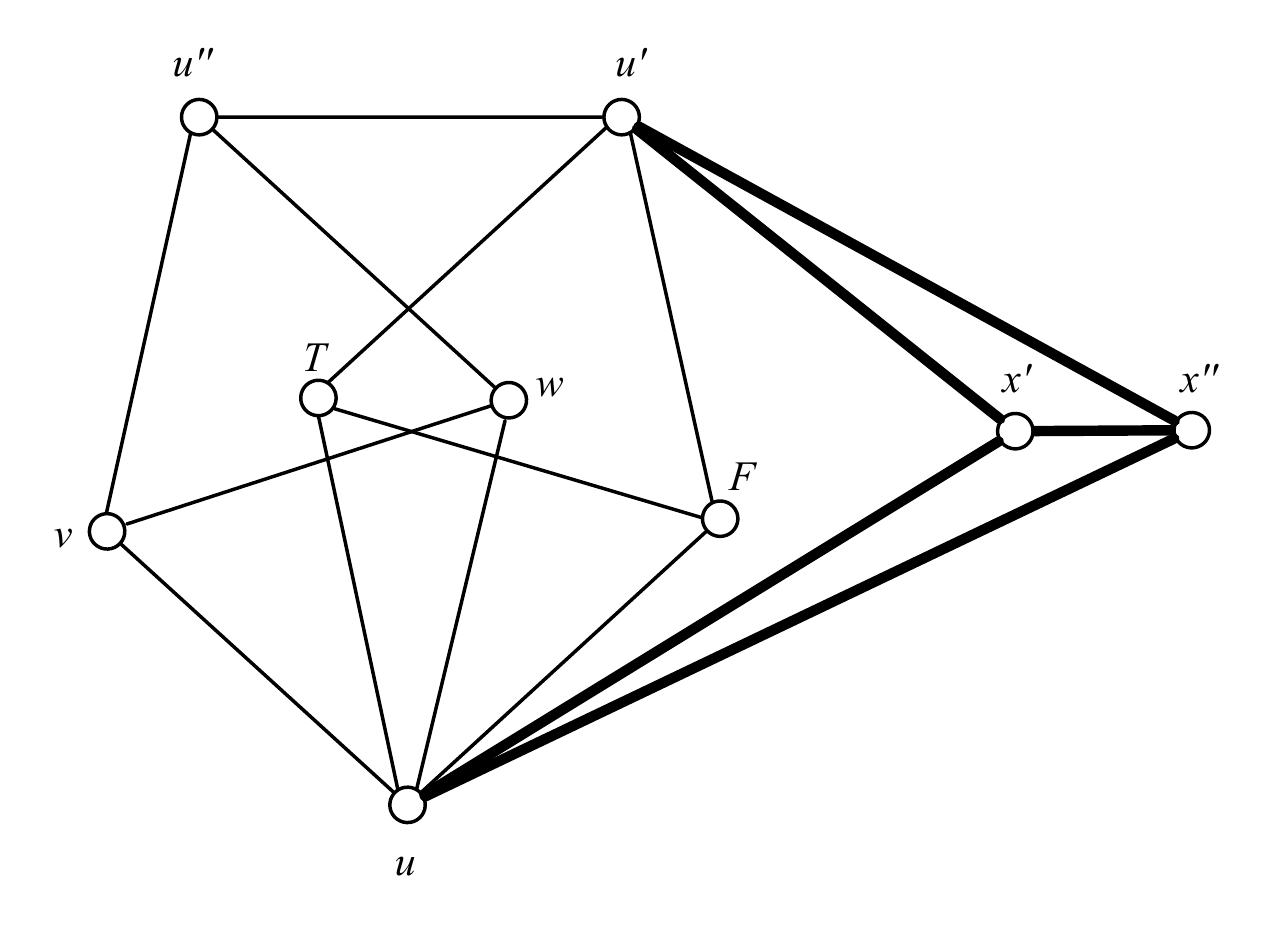,width=7cm}}
\caption{The variable part of the construction of $G_{\Phi}$.}
\label{fig:variable}
\end{figure}

The cornerstone of the construction is a copy of the Moser graph with vertices $\{u,T,F,u',u'',v,w\}$ and edges as drawn in the shaded part of the Fig.~\ref{fig:variable}. For every variable $x \in X$, we add an edge $x' \sim x''$, where $x'$ represents the positive literal and $x''$ the negation of $x$. Also each $x'$ and $x''$ are made adjacent to both $u$ and $u'$. The part of the construction of $G_{\Phi}$ that corresponds to a variable $x \in X$ can be seen in the bolded part of the Fig.~\ref{fig:variable}.

For every clause $c \in C$, say $c = (\lambda_{c,1} \vee \lambda_{c,2} \vee \lambda_{c,3})$, the clause gadget is depicted in Fig.~\ref{fig:clausegadget}. It contains six vertices, three of them being connected to the literals appearing in $c$ (one to each) and one adjacent to the special vertex $F$ of the Moser graph. 

Formally,
$$
V(G_{\Phi}) = \{u, v, w, u', u'', T, F\} \cup \bigcup_{x \in X} \{x', x''\} \cup \bigcup_{c \in C} \{c_{1}, c_{2}, c_{3}, c_{4}, c_{12}, c_{34}\}
$$
and 
\begin{eqnarray*}
E(G_{\Phi}) & = & \{u \sim T, u \sim F, T \sim F, T \sim u', F \sim u', u \sim v, u \sim w, v \sim w, \\
 &  & \indent  v \sim u'', w \sim u'', u' \sim u'' \}  \cup \\
 &  & \bigcup_{x \in X} \{x' \sim x'', x' \sim u, x'' \sim u, x' \sim u', x'' \sim u'\} \cup \\
 &  & \bigcup_{c \in C} \{c_{1} \sim c_{2}, c_{2} \sim c_{3}, c_{3} \sim c_{4}, c_{4} \sim c_{1}, c_{1} \sim c_{12}, c_{2} \sim c_{12}, c_{3} \sim c_{34},  \\
 &  & \indent  c_{4} \sim c_{34},  c_{12} \sim F, c_{2} \sim \lambda_{c,1}, c_{3} \sim \lambda_{c,2}, c_{34} \sim \lambda_{c,3}\}.
\end{eqnarray*}

\begin{figure}[ht]
\centerline{\epsfig{file=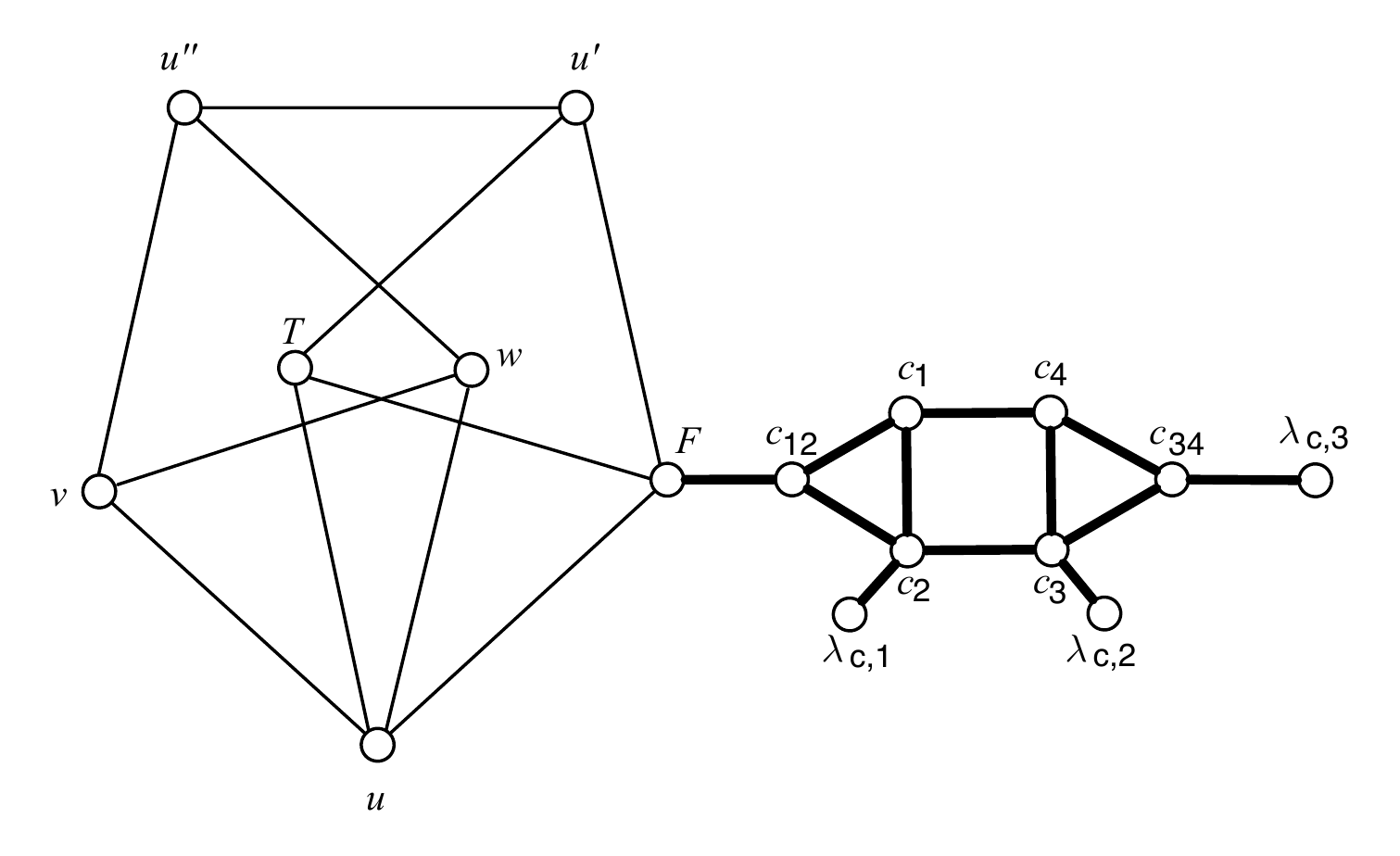,width=8cm}}
\caption{The clause part of $G_{\Phi}$.}
\label{fig:clausegadget}
\end{figure}

Suppose first that $G_{\Phi}$ is homomorphic to a graph $H$ with $\dim(H) = 2$ (to a planar unit distance graph $H$). Consider a planar unit distance realization of $H$. There are at most two paths of length two between two vertices of a planar unit distance graph \cite{UDGMAEHARA4}. Since there are four paths of length two between $u$ and $u'$ in $G_{\Phi}$, the homomorphism either places $x'$ in $T$ and $x''$ in $F$, or vice versa.
This placement defines a truth assignment on the variables of $\Phi$ - we say that $x$ is {\sf true} if $x'$ is placed in $T$, and that it is {\sf false} otherwise. We claim that $\Phi$ is satisfied by this assignment. Suppose there is an unsatisfied clause, say $c = (\lambda_{c,1} \vee \lambda_{c,2} \vee \lambda_{c,3})$. Then all three vertices $\lambda_{c,1}, \lambda_{c,2}, \lambda_{c,3}$ must be placed in $F$, and the clause gadget must map onto another copy of the Moser graph (with $F$ being its degree four vertex). But then, in the plane, the edge $c_2 \sim c_3$ cannot have a unit length. We have a contradiction.
% BH maybe another picture containing both copies of the Moser graph and the additional edge that prevents having a planar unit distance realization

\begin{figure}[ht]
\centerline{\epsfig{file=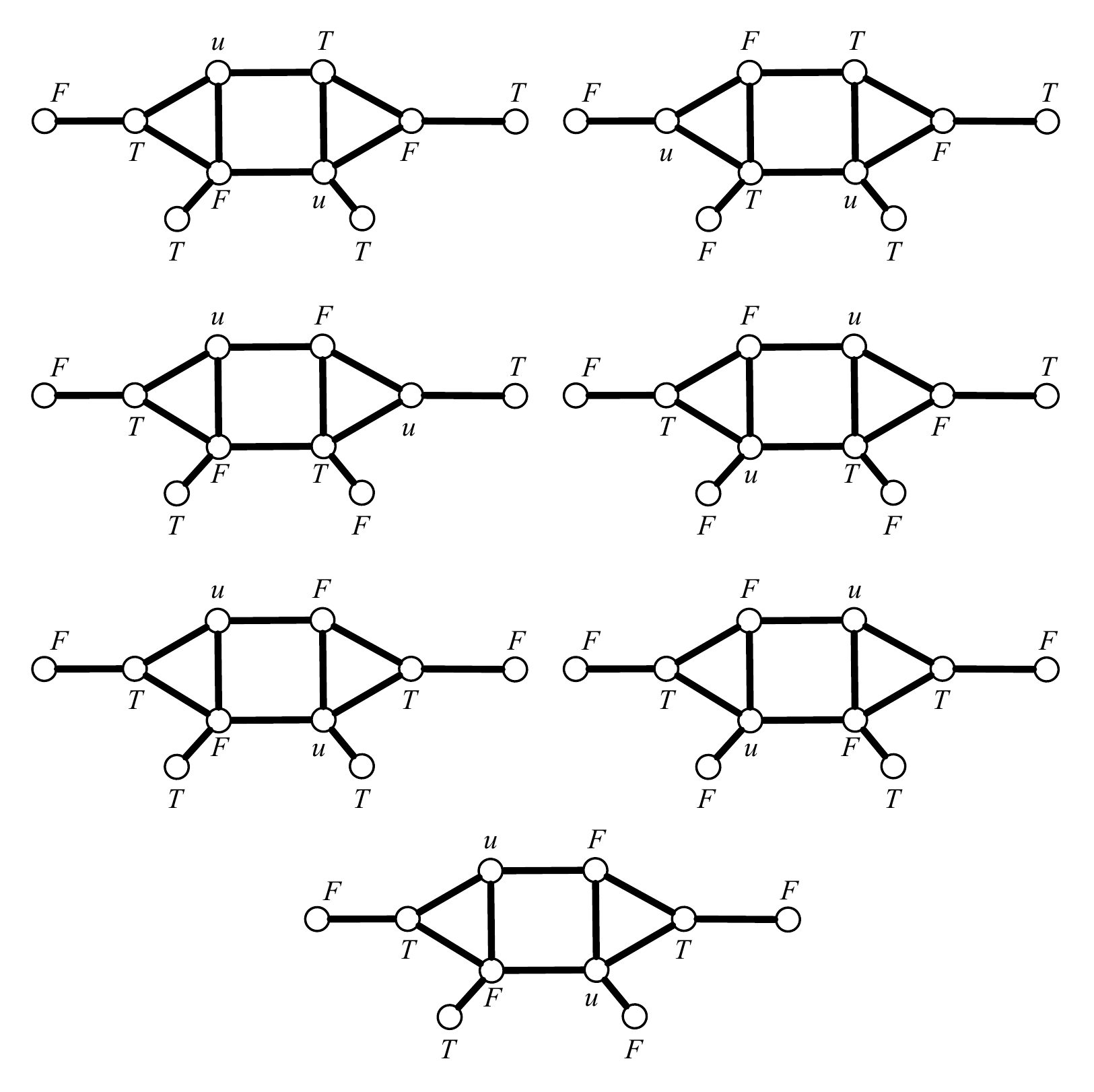,width=9cm}}
\caption{Case analysis of placement of the vertices of clause gadgets.}
\label{fig:placement}
\end{figure}

Suppose now that $\Phi$ is satisfied by a truth assignment $\phi: X\longrightarrow \{{\sf true, false}\}$. We fix a strict planar unit distance realization of the Moser graph, place $x'$ in $T$ and $x''$ in $F$ whenever $\phi(x)={\sf true}$ (and vice versa if $\phi(x)={\sf false}$). Following the case analysis in Fig.~\ref{fig:placement}, the vertices of the clause gadgets can be placed in vertices $u,T$ and $F$. Hence, we have constructed a well defined surjective homomorphism of $G_{\Phi}$ onto the Moser graph. Therefore, $G_{\Phi}$ is homomorphic to a graph with the Euclidean dimension 2.
\end{proof}

%%%%%%%%%%%%%%%%%%%%%%%%%%%%%%
%
%  $k$ dimensional unit distance representations
%
%%%%%%%%%%%%%%%%%%%%%%%%%%%%%%

\section{$k$ dimensional unit distance representations}

% Supporting lemma and proposition for Theorem 2
Let $\text{\pS}_{k}(\vec{s}, r)$ denote the $k$ dimensional (hiper)sphere in $\bR^{k}$ with center in $\vec{s}$ and radius $r$. When the center and the radius of a sphere are not important, the abbrevation $\text{\pS}_{k}$ will be used.

\begin{lemma}\label{le:presekksferjekm1}
Let $k>1$ be a natural number. A non-empty non-degenerated intersection of two $k$ dimensional spheres with distinct centers is $k-1$ dimensional sphere.
\end{lemma}
\begin{proof}
Let $c, r, R > 0$ be positive real numbers. We can, without loss of generality, assume, that the first sphere is centered at the origin $\vec{0}$ and the second one is centered on the $x_{1}$ axes in $\vec{c} = (c, 0, \ldots, 0)$.
Spheres $\text{\pS}_{k}(\vec{0},r)$ and $\text{\pS}_{k}(\vec{c},R)$ are determined by equations
$$x_1^2+x_2^2+ \ldots +x_k^2 = r^2$$
and
$$(x_1-c)^2+x_2^2+ \ldots +x_k^2 = R^2.$$
We subtract the above equations to obtain a linear equation for $x_1$ with a solution $a := {r^2+c^2-R^2 \over 2c}$. When $a$ is inserted into the first equation, the equality $x_2^2+ \ldots +x_k^2=r^2-a^2$ is obtained.
When $r^2-a^2 < 0$, the system does not have any solution and this case corresponds to a situation, where spheres do not intersect.
A situation $|a| = r$ corresponds to spheres that are tangent to each other and their intersection is a point - a degenerated sphere.
Assume now, that $r^2-a^2 > 0$. The equation $x_2^2+ \ldots +x_k^2 = r^2 - a^2$ determines the $k-1$ dimensional sphere in the hyperplane $x_1 = a$ in $\bR^{k}$.
\end{proof}

\begin{lemma}\label{lm:chileq3}
Let $\text{\pS}_{2}$ denote the circle in the Euclidean plane, which is circumscribed to a unit distance realization of the complete graph $K_{3}$ on three vertices. 
Let $G$ be a connected graph with (possibly degenerate) planar unit distance representation, which places all vertices of $G$ into points that lay on $\text{\pS}_{2}$. Then $\chi(G) \leq 3$.
\end{lemma}
\begin{proof}
Graph $K_{3}$ has unique (unique up to isometries) unit distance realization in the plane. It can be represented as an equilateral triangle with all sides of unit length, which has unique circumscribed circle. Denote the circumscribed circle by $\text{\pS}_{2}$. Every point on $\text{\pS}_{2}$ has exactly two points that are unit distance appart from it and belong to $\text{\pS}_{2}$. Graph $G$ has a unit distance representation on $\text{\pS}_{2}$ and can hence occupy only three distinct points on $\text{\pS}_{2}$, which are exactly the vertices of (possibly another) equilateral triangle with all sides of length one. Thus, there exists a proper 3 coloring of $G$ - colors being the vertices of the equlateral triangle.
\end{proof}

Lemma~\ref{lm:chileq3} can be generalized. Let $G_{k, r, \alpha}$ denote the graph with vertices being points of $k$ dimensional sphere $S_{k}(\vec{0}, r)$ with radius $r$ in $\bR^{k}$, where two vertices are connected if and only if they are at distance $\alpha$. L.~Lov\'asz \cite{Lovasz:1983p585} proved the following inequalities.

\begin{theorem}\label{th:chileqlovasz}
Let $k \geq 3$ be a natural number. Then, for $0 \leq \alpha \leq 2$ holds
$$
k \leq \chi(G_{k,1,\alpha})
$$
and
$$
\chi\left(G_{k,1,\sqrt{2(k+1) \over k}}\right) \leq k+1.
$$
\end{theorem}

\begin{corollary}\label{cor:chileqk}
Let $k \geq 3$ be a natural number and let $\text{\pS}_{k}$ denote the $k$ dimensional sphere in $\bR^{k}$, which is circumscribed to the unit distance realization of the complete graph $K_{k+1}$ on $k+1$ vertices. 
Let $G$ be a connected graph with (possibly degenerate) $k$ dimensional unit distance representation, which places all vertices of $G$ into points that lay on $\text{\pS}_{k}$. Then $\chi(G) \leq k+1$.
\end{corollary}
\begin{proof}
It is known, that the circumradius of the (hiper)sphere that is circumscribed to the regular simplex with $k+1$ vertices and all sides of length $\ell$, equals to $\ell \sqrt{k \over 2 (k+1)}$. 
Thus, the circumradius of sphere $\text{\pS}_{k}$ equals to $\sqrt{k \over 2 (k+1)}$. The representation of graph $G_{k, \sqrt{k \over 2 (k+1)}, 1}$ on sphere $\text{\pS}_{k}$ can be (down)scaled to obtain the representation on $k$ dimensional unit sphere (with radius one), that is circumscribed to the regular simplex with all sides of length $\sqrt{2(k+1) \over k}$.
Following Theorem~\ref{th:chileqlovasz}, $\chi\left(G_{k,1,\sqrt{2(k+1) \over k}}\right) \leq k+1$. The proper $k+1$ coloring of $G_{k,1,\sqrt{2(k+1) \over k}}$ gives rise to the proper $k+1$ coloring of $G_{k, \sqrt{k \over 2 (k+1)}, 1}$. Since $G$ is a subgraph of $G_{k, \sqrt{k \over 2 (k+1)}, 1}$, $\chi(G) \leq k+1$.
\end{proof}

% Theorem 2
\begin{theorem}
Let $k \geq 3$ be a natural number. Deciding if an input graph $G$ is homomorphic to a graph $H$ with $\dim(H) = k$ (with $\e(H) = k$) are NP-hard problems.
%Deciding if an input graph is \hsu\ (\hwu) in the 3-dimensional space are NP-hard problems.
\end{theorem}

\begin{proof}
Let $K'_{k}$ and $K''_{k}$ be two copies of the complete graph with $k \geq 3$ vertices. Let $v,w',w''$ be additional vertices, such that $v,w',w'' \notin V(K'_{k}) \cup V(K''_{k})$. 
Denote $M'_{k} = K'_{k} \join \{v, w'\}$, $M''_{k} = K''_{k} \join \{v, w''\}$ and $M_{k} = M'_{k} \union M''_{k} \union \{w' \sim w''\}$, e.g. Figure~\ref{fig:kcoldeg1}. Graph $M_{k}$ is well known Moser-Raiskii spindle.

\begin{figure}[ht]
\centerline{\epsfig{file=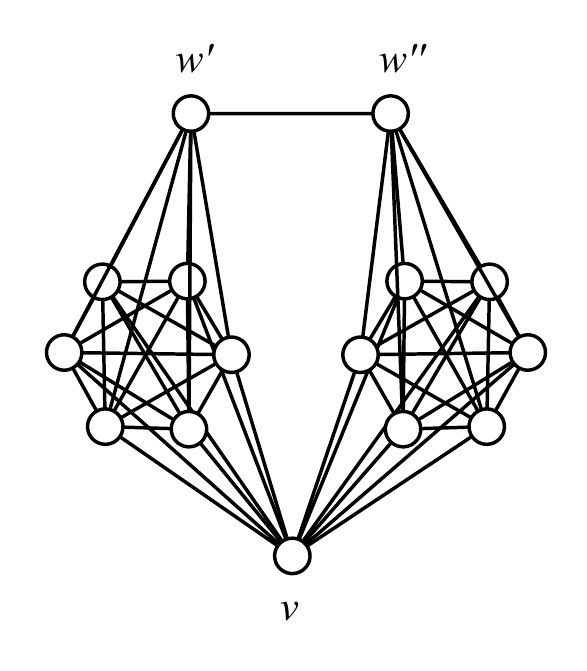,width=4.5cm}}
\caption{The rigid graph $M_{6}$ contains two similar subgraphs, which are obtained as a graph join of the complete graph on six vertices and two disconnected vertices.}
\label{fig:kcoldeg1}
\end{figure}

Note that $\e(K'_{k}) = \dim(K'_{k}) = k-1$. Since $K'_{k} \join v$ is the complete graph on $k+1$ vertices, $\dim(K'_{k} \join v) = \e(K'_{k} \join v) = k$. Applying another graph join does not change any of dimensions. Thus,  $e(M'_{k}) = e(M''_{k}) = \dim(M'_{k}) = \dim(M''_{k}) = k$. Both $M'_{k}$ and $M''_{k}$ are rigid  in $\bR^{k}$. % rigid in a sense that is not flexible
Consider a unit distance realization of $M'_{k} \union M''_{k}$ in $\bR^{k}$. We can rotate $M''_{k}$ around the representation of the vertex $v$, such that $w'$ and $w''$ results distance one appart.
Hence $\dim(M_{k}) = k$. Additionally, we can rotate the subgraph $K''_{k}$ in such a way, that $\e(M_{k})=k$. 
Note that $H_{k}$ has the property that every $k$ dimensional unit distance realization in $\bR^{k}$, places vertices $v$ and $w'$ in distinct points (this is guaranteed by the edge $w' \sim w''$). 

Let us construct the graph $H_{k}$, such that $V(H_{k}) := V(G) \cup V(M_{k})$ and 
$$
E(H_{k}) := E(G) \union E(M_{k}) \union \bigcup_{u \in V(G)} \{u \sim v, u \sim w'\},
$$
where vertices $v, w'$ are vertices from $M_{k}$.

We reduce from graph \COL. We will prove that, when $H_{k}$ is homomorphic to a graph with the dimension $k$,
then $\chi(G) \leq k$, while from $\chi(G) \leq k$ follows, that $H_{k}$ is homomorphic to a graph with the Euclidean dimension $k$.
% e leq k=> dim leq k => chi leq k 
% chi leq k => e leq k => dim leq k

Suppose $H_k$ is homomorphic to a graph with the dimension $k$. Then $H_k$ has (possibly degenerate) unit distance representation $\rho$ in $\bR^{k}$. All vertices of $G$ are one appart from $\rho(v)$ and from $\rho(w')$, and they all belong to the intersection of two unit $k$ dimensional spheres with centers in $\rho(v)$ and $\rho(w')$. Using Lemma~\ref{le:presekksferjekm1}, a non-empty non-degenerated intersection of two $k$ dimensional spheres with distinct centers is $k-1$ dimensional sphere. All vertices of the complete graph $K'_{k}$ belong to the same intersection and hence to the same $k-1$ dimensional sphere. 
If $G$ is connected, we can use Corollary~\ref{cor:chileqk} and $\chi(G) \leq k$. If $G$ is not connected, we can observe its components in similar matter.

Assume now, that $\chi(G) \leq k$. Using proper $k$-coloring of $G$ we can map vertices of $G$ into vertices of $K'_{k}$, see \cite{GAHHELL}.
Hence a homomorphic image of $H_{k}$ is $M_{k}$ and $\e(M_{k}) = k$.

As any $k$ dimensional strict unit distance graph is a $k$ dimensional unit distance graph, it follows that for a given graph $G$, $\dim(G) \leq \e(G)$. Hence we proved both claims of this Theorem.
\end{proof}

%%%%%%%%%%%%%%%%%%%%%%%%%%%%%%
%
%  Concluding remarks
%
%%%%%%%%%%%%%%%%%%%%%%%%%%%%%%

\section{Concluding remarks}
There are several results in connection to computational complexity and unit distance graphs.
Firstly, it is known, that a graph reconstruction problem is strongly NP-complete in $1$ dimension and strongly NP-hard in higher dimensions, see \cite{saxe}. For general edge-lenghts, graph realization problems and problems concerning rigidity of structures in connection to computational complexity, were observed in \cite{Hendrickson:1992p3225} and \cite{Yemini:1979p3248}.

P.~Eades and S.~Whitesides proved in \cite{Eades:1996p3245} that the decision problem: ``Is a planar graph $G$ a matchstick graph?'' is NP-hard. They used a reduction from the well known NP-complete problem \NAE3SAT, a special instance of the \SAT problem, in which, each term has to have at least one true and at least one false literal. Their proof was constructive, they used a matchstick graph to construct a mechanical device, called {\em logic engine}, that was used to simulate the \NAE3SAT problem instance.
Cabello {\em et al.} in \cite{Cabello:2004p560} used similar approach to prove that the decision problem: ``Is a planar $3$-connected infinitesimally rigid graph $G$ a matchstick graph?'' is NP-hard. 

Additionally, P.~Eades and N.~C.~Wormald proved in \cite{Eades:1990p3289} that the decision problem: ``Is a $2$-connected graph $G$ a matchstick graph?'' is NP-hard, too. 

% BH This section is not finished. Maybe we should add some future work?

% ---------------------------------------------------------------
% ------------------======== References ========-----------------
% ---------------------------------------------------------------


\begin{thebibliography}{99}

\bibitem{BobenPisanski}
M.~Boben, T.~Pisanski.
Polycyclic configurations.
\emph{European J. Combin.} \textbf{24} (2003), 4:431--457.

\bibitem{BuckleyHarary}
F.~Buckley, F.~Harary.
On the Euclidean dimension of a Wheel.
\emph{Graphs Combin.} \textbf{4} (1988), 23--30.

\bibitem{Cabello:2004p560}
S.~Cabello, E.~Demaine, G.~Rote.
Planar embeddings of graphs with specified edge lengths.
\emph{Lect. Notes Comput. Sc.} \textbf{2912} (2004), 283--294.

\bibitem{Eades:1996p3245}
P.~Eades, S.~Whitesides.
The logic engine and the realization problem for nearest neighbor graphs.
\emph{Theor. Comput. Sc.} \textbf{169} (1996), 1:23--37.

\bibitem{Eades:1990p3289}
P.~Eades, N.~C.~Wormald.
Fixed edge-length graph drawing is NP-hard.
\emph{Discrete Appl. Math.} \textbf{28} (1990), 2:111--134.

\bibitem{ErdosHararyTutte}
P.~Erd\"os, F.~Harary, W.~T.~Tutte.
On the dimension of a graph.
\emph{Mathematika} \textbf{12} (1965), 118--122.

\bibitem{pearlsgt}
N.~Hartsfield, G.~Ringel.
\emph{Pearls in Graph Theory: A Comprehensive Introduction}. Revised and augmented.
Academic Press, San Diego, 1994.

\bibitem{GAHHELL}
P.~Hell, J.~Ne{\v{s}}et{\v{r}}il.
\emph{Graphs and Homomorphisms}. 
Oxford University Press, 2004.

\bibitem{Hendrickson:1992p3225}
B.~Hendrickson.
Conditions For Unique Graph Realizations.
\emph{SIAM J. Comput.} \textbf{21} (1992), 1:65--84.

\bibitem{HorvatPisanski}
B.~Horvat, T.~Pisanski.
Unit distance representations of the Petersen graph in the plane.
\emph{Ars Combin.} (to appear).

\bibitem{HorvatPisanski2}
B.~Horvat, T.~Pisanski.
Products of unit distance graphs.
\emph{Discrete Math.} (to appear).

\bibitem{Lovasz:1983p585}
L.~Lov\'asz.
Self-dual polytopes and the chromatic number of distance graphs on the sphere.
\emph{Acta Sci. Math. (Szeged)} \textbf{45} (1983), 1-4:317--323.

\bibitem{UDGMAEHARA3}
H.~Maehara.
Note on Induced Subgraphs of the Unit Distance Graph.
\emph{Discrete Comput. Geom.} \textbf{4} (1989), 15--18.

\bibitem{UDGMAEHARA4}
H.~Maehara.
On the Euclidean dimension of a complete multipartite graph.
\emph{Discrete Math.} \textbf{72} (1988), 285--289.

\bibitem{UDGMAEHARA2}
H.~Maehara, V.~R\"{o}dl.
On the Dimension to Represent a Graph by a Unit Distance Graph.
\emph{Graphs Combin.} \textbf{6} (1990), 365--367.

\bibitem{saxe}
J.~B.~Saxe.
Embeddability of weighted graphs in k-space is strongly NP-hard.
In \emph{Proc. 17th Al lerton Conf. Commun. Control Comput.}, 480--489, 1979. 

\bibitem{Yemini:1979p3248}
Y.~Yemini.
Some theoretical aspects of position-location problems.
In \emph{Proc. 20th Annu. IEEE Sympos. Found. Comput. Sci.}, 1--8, 1979. 


\end{thebibliography}
\end{document}